\documentstyle[12pt]{article}
\begin{document}
\setlength{\textheight}{574pt}
\setlength{\textwidth}{432pt}
\setlength{\oddsidemargin}{18pt}
\setlength{\topmargin}{14pt}
\setlength{\evensidemargin}{18pt}
\newtheorem{theorem}{Theorem}[section]
\newtheorem{lemma}{Lemma}[section]
\newtheorem{corollary}{Corollary}[section]
\newtheorem{remark}{Remark}[section]
\newtheorem{definition}{Definition}[section]
\newtheorem{problem}{Problem}
\newtheorem{proposition}{Proposition}[section]
\title{{\bf BOUND OF AUTOMORPHISMS OF PROJECTIVE VARIETIES OF GENERAL TYPE}}
\date{April, 2000}
\author{Hajime TSUJI}
\maketitle
\begin{abstract}
We prove that there exists a positive integer $C_{n}$ depending
only on $n$ such that for every smooth projective $n$-fold of 
general type $X$ defined over {\bf C}, the  automorphism group 
$\mbox{Aut}(X)$ of $X$ satisfies 
\[
\sharp\mbox{Aut}(X) \leq C_{n}\cdot\mu (X,K_{X}),
\]
where $\mu (X,K_{X})$ is the volume of $X$ with respect to 
$K_{X}$.  MSC14E05,32J25. 
\end{abstract}
\section{Introduction}
The automorphim group of a projective variety of general type 
is known to be finite. 
For every  curve $C$ of genus $g \geq 2$, we have 
the estimate : 
\[
\sharp\mbox{Aut}(C) \leq 84(g-1)
\]
by well known Hurwitz's theorem. 

In the case of surfaces, G. Xiao proved that 
for every smooth minimal surface of general type 
\[
\sharp\mbox{Aut}(S) \leq 1764\cdot K_{S}^{2}
\]
holds \cite{x1,x2}.
The main purpose of this article is to prove the following theorem.

\begin{theorem} There exists a positive number $C_{n}$ which depends
only on $n$ such that for every smooth projective $n$-fold $X$
of general type defined over complex numbers, the automorphism group 
$\mbox{Aut}(X)$  of $X$ satisfies the estimate:
\[
\sharp\mbox{Aut}(X) \leq C_{n}\cdot\mu (X,K_{X}),
\]
where $\mu (X,K_{X})$ is the volume of $X$ with respect to 
$K_{X}$ (cf. Definition 2.3).
\end{theorem}
The method of the proof of Theorem 1.1 is a combination of 
the ideas in  \cite{x1,x2} and \cite{t4}.  
Let $X$ be a projective $n$-fold of general type 
and let $G$ denotes the automorphism group of $X$. 
Since $G$ acts on the canonical ring $R(X,K_{X})$ 
of $X$, 
by \cite{tu3} we may assume that $X$ is a canonical 
model, i.e. $X$ has only canonical singularity and 
$K_{X}$ is ample (our proofs of Theorem 1.1 and Therem 1.2 below depend 
on the finite generation of canonical rings of varieties of general type in 
\cite{tu3} which has not yet published. 
For the safe side, one may restrict oneselves to the case of $\dim X \leq 3$ 
(cf. \cite{mo}))
The quotient $X/G$ is a projective variety. 
Let $K_{X/G,orb}$ be the orbifold canonical divisor 
of $X/G$. 
Then we see that 
\[
\mid G\mid = K_{X}^{n}/K_{X/G,orb}^{n}
\]
holds, where $\mid G\mid$ denotes the order of $G$. 
Since $\mu (X,K_{X}) = K_{X}^{n}$ holds in this case,
we see that Theorem 1.1 follows from the following 
theorem. 
\begin{theorem} 
Let $X$, $G$ be as above. 
There exists a positive constant $c_{n}$ depending only 
on $n$ such that 
\[
K_{X/G,orb}^{n} \geq c_{n}
\]
holds.
\end{theorem}
It is easy to see $c_{1}$ can be taken 
to be $1/42$. 
This leads to Hurwicz's theorem. 
G. Xiao proved that 
$c_{2}$ can be taken as $1/1764$ (\cite{x1,x2}).

The key ingredient of the proof of Theorem 1.2 is 
the subadjunction formula in \cite{ka} which relates 
the canonical divisor of the minimal center of logcanonical 
singularities and the canonical divisor of the ambient space.   
Using this we  see that $X/G$ with 
$\mu (X/G,K_{X/G,orb}) = K_{X/G,orb}^{n} \leq 1$ is
birationally bounded by the inductive procedure in 
\cite{t4}. 
Then Theorem 1.1 and Theorem 1.2 follows from a Diophantine consideration. 

Therem 1.1 and Theorem 1.2 are not effective in the sense that 
there exist no explicit estimates of $C_{n}$ and $c_{n}$.

\section{Preliminaries}
\subsection{Orbifold canonical divisors}
Let $X$ be a projective variety of general type
with only canonical singularities. 
Let $G$ denote the automorphism group of $X$. 
Iis well known that $G$ is a finite group. 
The quotient $X/G$ is a projective variety.
Let $\tilde{X}$ be the equivalent resolution of $X$ 
with respect to $G$ such that 
$\tilde{X}/G$ is also smooth.
We may take $\tilde{X}$ such that the ramification divisor $R$   
of 
\[
\tilde{\pi} : \tilde{X} \longrightarrow \tilde{X}/G
\]
and the branch locus $B = (\tilde{\pi}_{*}(R))_{red}$ is a divisor with normal crossings. 
Let $B = \sum_{i} B_{i}$ be the irreducible decomposition of 
$B$. 
Then there exists a set of positive integers $m_{i}$ such that 
\[
K_{\tilde{X}} = \tilde{\pi}^{*}(K_{\tilde{X}/G} +\sum_{i}\frac{m_{i}-1}{m_{i}}B_{i}) 
 \]
Let 
\[
\varpi : \tilde{X}/G \longrightarrow X/G
\]
be the natural morphism. 
We set 
\[
K_{X/G,orb} := \varpi_{*}(K_{\tilde{X}/G} +\sum_{i}\frac{m_{i}-1}{m_{i}}B_{i})
\]
and call it the obifold canonical divisor of $X/G$. 
Let 
\[
\pi : X \longrightarrow X/G
\]
be the natural morphism.  
Then 
\[
K_{X} = \pi^{*}K_{X/G,orb} 
\] 
holds. 
The orbifold canonical ring is defined by 
\[
R(X/G,K_{X/G,orb}):= R(X,K_{X})^{G}.
\]
And the linear system $\mid [mK_{X/G,orb}]\mid$ is given by 
\[
\mid [mK_{X/G,orb}]\mid = \mid mK_{X}\mid^{G}.
\]
Hence we have that 
\[
R(X/G,K_{X/G,orb}) = \oplus_{m\geq 0}\Gamma (X/G,{\cal O}_{X/G}([mK_{X/G,orb}]))\]
holds. 

\subsection{Multiplier ideal sheaves}
In this section, we shall review the basic definitions and properties
of multiplier ideal sheaves.
\begin{definition}
Let $L$ be a line bundle on a complex manifold $M$.
A singular hermitian metric $h$ is given by
\[
h = e^{-\varphi}\cdot h_{0},
\]
where $h_{0}$ is a $C^{\infty}$-hermitian metric on $L$ and 
$\varphi\in L^{1}_{loc}(M)$ is an arbitrary function on $M$.
\end{definition}
The curvature current $\Theta_{h}$ of the singular hermitian line
bundle $(L,h)$ is defined by
\[
\Theta_{h} := \Theta_{h_{0}} + \sqrt{-1}\partial\bar{\partial}\varphi ,
\]
where $\partial\bar{\partial}$ is taken in the sense of a current.
The $L^{2}$-sheaf ${\cal L}^{2}(L,h)$ of the singular hermitian
line bundle $(L,h)$ is defined by
\[
{\cal L}^{2}(L,h) := \{ \sigma\in\Gamma (U,{\cal O}_{M}(L))\mid 
\, h(\sigma ,\sigma )\in L^{1}_{loc}(U)\} ,
\]
where $U$ runs opens subsets of $M$.
In this case there exists an ideal sheaf ${\cal I}(h)$ such that
\[
{\cal L}^{2}(L,h) = {\cal O}_{M}(L)\otimes {\cal I}(h)
\]
holds.  We call ${\cal I}(h)$ the multiplier ideal sheaf of $(L,h)$.
If we write $h$ as 
\[
h = e^{-\varphi}\cdot h_{0},
\]
where $h_{0}$ is a $C^{\infty}$ hermitian metric on $L$ and 
$\varphi\in L^{1}_{loc}(M)$ is the weight function, we see that
\[
{\cal I}(h) = {\cal L}^{2}({\cal O}_{M},e^{-\varphi})
\]
holds.
We have the following vanishing theorem.

\begin{theorem}(Nadel's vanishing theorem \cite[p.561]{n})
Let $(L,h)$ be a singular hermitian line bundle on a compact K\"{a}hler
manifold $M$ and let $\omega$ be a K\"{a}hler form on $M$.
Suppose that $\Theta_{h}$ is strictly positive, i.e., there exists
a positive constant $\varepsilon$ such that
\[
\Theta_{h} \geq \varepsilon\omega
\]
holds.
Then ${\cal I}(h)$ is a coherent sheaf of ${\cal O}_{M}$-ideal 
and for every $q\geq 1$
\[
H^{q}(M,{\cal O}_{M}(K_{M}+L)\otimes{\cal I}(h)) = 0
\]
holds.
\end{theorem}

\subsection{Analytic Zariski decomposition}

To study a big line bundle we introduce the notion of analytic Zariski
decompositions.
By using analytic Zariski decompositions, we can handle big line bundles
like a nef and big line bundles.
\begin{definition}
Let $M$ be a compact complex manifold and let $L$ be a line bundle
on $M$.  A singular hermitian metric $h$ on $L$ is said to be 
an analytic Zariski decomposition, if the followings hold.
\begin{enumerate}
\item $\Theta_{h}$ is a closed positive current,
\item for every $m\geq 0$, the natural inclusion
\[
H^{0}(M,{\cal O}_{M}(mL)\otimes{\cal I}(h^{m}))\rightarrow
H^{0}(M,{\cal O}_{M}(mL))
\]
is isomorphim.
\end{enumerate}
\end{definition}
\begin{remark} If an AZD exists on a line bundle $L$ on a smooth projective
variety $M$, $L$ is pseudoeffective by the condition 1 above.
\end{remark}

\begin{theorem}(\cite{tu,tu2})
 Let $L$ be a big line  bundle on a smooth projective variety
$M$.  Then $L$ has an AZD. 
\end{theorem}

\subsection{Volume of projective varieties}

To measure the positivity of big line bundles 
on a projective variety we shall introduce a volume of a projective variety
with respect to a line bundle.

\begin{definition} Let $L$ be a line bundle on a compact complex 
manifold $M$ of dimension $n$. 
We define the $L$-volume of $M$ by
\[
\mu (M,L) := n!\cdot\limsup_{m\rightarrow\infty}m^{-n}
\dim H^{0}(M,{\cal O}_{M}(mL)).
\]
\end{definition}
\begin{definition}(\cite{tu3})
Let $L$ be a big line bundle on a smooth projective variety
$X$. 
Let $Y$ be a subvariety of $X$ of dimension $r$.
We define the volume $\mu (Y,L)$ of $Y$ with respect to 
$L$ by 
\[
\mu (Y,L) := r!\cdot\limsup_{m\rightarrow\infty}m^{-r}
\dim H^{0}(Y,{\cal O}_{Y}(mL)\otimes{\cal I}(h^{m})/tor),
\]
where $h$ is an AZD of $L$ and $tor$ denotes the torsion part of the sheaf ${\cal O}_{Y}(mL)\otimes{\cal I}(h^{m})$.
This definition can be easily generalized to the case that 
$L$ is a {\bf Q}-line bundle.  
\end{definition}

\section{Stratification of varieties by multiplier ideal sheaves}

Let $X$ be a smooth projective $n$-fold of general type.
Then the canonical ring $R(X,K_{X})$ is finitely generated by \cite{tu3}.  
Let $X_{can}$ be the canonical model of $X$. 
$K_{X_{can}}$ is an ample ${\bf Q}$-Cartier divisor on $X_{can}$.
We assume that the natural rational map 
\[
\varphi : X -\cdots \rightarrow X_{can} 
\]
is a morphism.
Let $h_{can}$ be a $C^{\infty}$-hermitian metric 
on $K_{X_{can}}$ induced from the Fubini-Study metric on the hyperplane bundle of a projective space by a projective embedding of $X_{can}$ associated 
with $\mid rK_{X_{can}}\mid$ where  $r$ is a sufficiently large positive integer such that $rK_{X_{can}}$ is Cartier. 
Then $h_{can}$ has strictly positive curvature on $X_{can}$. 
$h_{can}$ induces a singular hermitian metric $h$ on $K_{X}$ in a natural manner.  
By the definition, $h$ is an AZD of $K_{X}$. 
To prove Theorem 1.1, we may replace $X$ by any birational model of $X$, we may assume that there exists an effective {\bf Q}-divisor $N$ such that 
\[
{\cal I}(h^{m}) = {\cal O}_{X}(-[mN])
\]
holds for every $m\geq 0$.
In particular we may and do assume that ${\cal I}(h^{m})$ is locally free for every $m\geq 0$. 
Let us denote $\mu (X/G,K_{X/G,orb})$ by $\mu_{0}$.
We set 
\[
X^{\circ} = \{ x\in X\mid \mbox{$\varphi$ is a local isomorphism around $x$}\} .
\]

Let $G$ be the group of the birational automorphism of 
$X$. 
To prove Theorem 1.1, we may assume that 
$G$ acts $X$ regularly and 
$X/G$ is also smooth. 
Let 
\[
\pi : X\longrightarrow X/G
\]
be the natural morphism. 
We set 
\[
(X/G)^{\circ} = \pi (X^{\circ}).
\]

\begin{lemma} Let $x,y$ be distinct points on $(X/G)^{\circ}$.  
We set 
\[
{\cal M}_{x,y} = {\cal M}_{x}\otimes{\cal M}_{y}
\]

Let $\varepsilon$ be a sufficiently small positive number.
Then 
\[
H^{0}(X/G,{\cal O}_{X/G}(mK_{X/G,orb})\otimes{\cal M}_{x,y}^{\lceil\sqrt[n]{\mu_{0}}
(1-\varepsilon )\frac{m}{\sqrt[n]{2}}\rceil})\neq 0
\]
for every sufficiently large $m$, where ${\cal M}_{x},{\cal M}_{y}$ denote the
maximal ideal sheaf of the points $x,y$ respectively.
\end{lemma}
{\em Proof of Lemma 3.1}.   
Let us consider the exact sequence:
\[
0\rightarrow H^{0}(X/G,{\cal O}_{X/G}(mK_{X/G,orb})\otimes
{\cal M}_{x,y}^{\lceil\sqrt[n]{\mu_{0}}(1-\varepsilon )\frac{m}{\sqrt[n]{2}}\rceil})
\rightarrow H^{0}(X/G,{\cal O}_{X/G}(mK_{X/G,orb}))\rightarrow
\]
\[
  H^{0}(X/G,{\cal O}_{X/G}
(mK_{X/G,orb})/{\cal M}_{x,y}^{\lceil\sqrt[n]{\mu_{0}}(1-\varepsilon )\frac{m}{\sqrt[n]{2}}\rceil}).
\]
Since 
\[
n!\limsup_{m\rightarrow\infty}m^{-n}\dim H^{0}(X/G,{\cal O}_{X/G}(mK_{X/G,orb})
/{\cal M}_{x,y}^{\lceil\sqrt[n]{\mu_{0}}(1-\varepsilon )\frac{m}{\sqrt[n]{2}}\rceil})
=
\mu_{0}(1-\varepsilon )^{n} < \mu_{0}
\]
hold, we see that Lemma 3.1 holds.  {\bf Q.E.D.}
\vspace{5mm}

Let us take a sufficiently large positive integer $m_{0}$ and let $\sigma$
be a general (nonzero) element of  
$H^{0}(X/G,{\cal O}_{X/G}(m_{0}K_{X/G,orb})\otimes
{\cal M}_{x,y}^{\lceil\sqrt[n]{\mu_{0}}(1-\varepsilon )\frac{m_{0}}{\sqrt[n]{2}}\rceil})$.
We define a singular hermitian metric $h_{0}$ on $K_{X/G,orb}$ by
\[
h_{0}(\tau ,\tau ) := \frac{\mid \tau\mid^{2}}{\mid \sigma\mid^{2/m_{0}}}.
\]
Then 
\[
\Theta_{h_{0}} = \frac{2\pi}{m_{0}}(\sigma )
\]
holds, where $(\sigma )$ denotes the closed positive current 
defined by the divisor $(\sigma )$.
Hence $\Theta_{h_{0}}$ is a closed positive current.
Let $\alpha$ be a positive number and let ${\cal I}(\alpha )$ denote
the multiplier ideal sheaf of $h_{0}^{\alpha}$, i.e.,
\[
{\cal I}(\alpha ) = 
{\cal L}^{2}({\cal O}_{X/G},(\frac{h_{0}}{h_{X/G}})^{\alpha}),
\]
where $h_{X/G}$ is an arbitrary $C^{\infty}$-hermitian metric on
$K_{X/G,orb}$.
Let us define a positive number $\alpha_{0} (= \alpha_{0}(x,y))$ by
\[
\alpha_{0} := \inf\{\alpha > 0\mid ({\cal O}_{X/G}/{\cal I}(\alpha ))_{x}\neq 0\,\mbox{and}\, ({\cal O}_{X/G}/{\cal I}(\alpha))_{y}\neq 0\}.
\]
Since $(\sum_{i=1}^{n}\mid z_{i}\mid^{2})^{-n}$ is not locally integrable 
around $O\in \mbox{{\bf C}}^{n}$, by the construction of $h_{0}$, we see
that 
\[
\alpha_{0}\leq \frac{n\sqrt[n]{2}}{\sqrt[n]{\mu_{0}}(1-\varepsilon )}
\]
holds.
Then one of the following two cases occurs. \vspace{10mm} \\
{\bf Case} 1.1:  For every small positive number $\delta$, 
${\cal O}_{X/G}/{\cal I}(\alpha_{0}-\delta )$  has $0$-stalk 
at both $x$ and $y$. \\
{\bf Case} 1.2: For every small positive number $\delta$, 
${\cal O}_{X/G}/{\cal I}(\alpha_{0}-\delta )$  has nonzero-stalk 
at one of $x$ or $y$ say $y$. \vspace{5mm} \\
First we consider Case 1.1.
Let $\delta$ be a sufficiently small positive number and  
let $V_{1}$ be the germ of subscheme at $x$ defined by the ideal sheaf 
${\cal I}(\alpha_{0}+\delta )$.
By the coherence of ${\cal I}(\alpha ) (\alpha > 0)$, we see that 
if we take $\delta$ sufficiently small, then $V_{1}$ is independent
of $\delta$.  It is also easy to verify that $V_{1}$ is reduced if 
we take $\delta$ sufficiently small. 
In fact if we take a log resolution of 
$(X/G,\frac{\alpha_{0}}{m_{0}}(\sigma ))$, 
$V_{1}$ is the image of the divisor with discrepancy $-1$ 
(for example cf. \cite[p.207]{he}). 
Let $(X/G)_{1}$ be a subvariety of $X/G$ which defines a branch of $V_{1}$
at $x$. 
We consider the following two cases. \vspace{10mm} \\
{\bf Case} 2.1: $(X/G)_{1}$ passes through both $x$ and $y$, \\
{\bf Case} 2.2: Otherwise \vspace{10mm} \\

For the first we consider Case 2.1.
Suppose that $(X/G)_{1}$ is not isolated at $x$.  Let $n_{1}$ denote
the dimension of $(X/G)_{1}$.  Let us define the volume $\mu_{1}$ of $(X/G)_{1}$
with respect to $K_{X/G,orb}$ by
\[
\mu_{1} := \mu ((X/G)_{1},K_{X/G,orb}).
\]
Since $x\in X/G^{\circ}$, we see that $\mu_{1} > 0$ holds.

\begin{lemma} Let $\varepsilon$ be a sufficiently small positive number and let $x_{1},x_{2}$ be distinct regular points on $(X/G)_{1}\cap X/G^{\circ}$. 
Then for a sufficiently large $m >1$ divisible by $\mid G\mid$,
\[
H^{0}((X/G)_{1},{\cal O}_{(X/G)_{1}}(mK_{X/G,orb})\otimes{\cal I}(h^{m})\otimes
{\cal M}_{x_{1},x_{2}}^{\lceil\sqrt[n_{1}]{\mu_{1}}(1-\varepsilon )\frac{m}{\sqrt[n_{1}]{2}}\rceil})\neq 0
\]
holds.
\end{lemma}
The proof of Lemma 3.2 is identical as that of Lemma 3.1, since 
\[
{\cal I}(h^{m})_{x_{i}} = {\cal O}_{X/G,x_{i}} (i=1,2)
\]
hold for every $m$.

By Kodaira's lemma there is an effective ${\bf Q}$-divisor $E$ such
that $K_{X/G,orb}- E$ is ample.
Let $\ell$ be a sufficiently large positive integer such that
\[
L := \ell (K_{X/G,orb}- E)
\]
is a line bundle and statisfies the property in Lemma 3.3.
\begin{lemma}
If we take $\ell$ sufficiently large, then 
\[
\phi_{m} : H^{0}(X/G,{\cal O}_{X/G}(mK_{X/G,orb}+L)\otimes{\cal I}(h^{m}))\rightarrow \]
\[
\hspace{40mm}
H^{0}((X/G)_{1},{\cal O}_{(X/G)_{1}}(mK_{X/G,orb}+L )\otimes{\cal I}(h^{m}))
\]
is surjective for every  $m\geq 0$ divisible by $\mid G\mid$.
\end{lemma}
{\em  Proof}.
Let us take a locally free resolution of the ideal sheaf ${\cal I}_{(X/G)_{1}}$
of $(X/G)_{1}$.
\[
0\leftarrow {\cal I}_{(X/G)_{1}}\leftarrow {\cal E}_{1}\leftarrow {\cal E}_{2}
\leftarrow \cdots \leftarrow {\cal E}_{k}\leftarrow 0.
\]
Then by the trivial extention of  
the case of vector bundles, if $r$ is sufficiently large, we see that
\[
H^{q}(X/G,{\cal O}_{X/G}(mK_{X/G,orb}+L)\otimes{\cal I}(h^{m})
\otimes{\cal E}_{j}) = 0
\]
holds for every $m\geq 1$, $q\geq 1$ and  $1\leq j\leq k$.
In fact if we take $\ell$ sufficiently large, we see that for every $j$, 
${\cal O}_{X/G}(L - K_{X/G})\otimes {\cal E}_{j}$ admits a $C^{\infty}$-hermitian metric $g_{j}$ such that
\[
\Theta_{g_{j}} \geq \mbox{Id}_{E_{j}}\otimes \omega
\]
holds, where $\omega$ is a K\"{a}hler form on $X/G$.
By \cite[Theorem 4.1.2 and Lemma 4.2.2]{ca} we have the desired vanishing. 

Hence 
\[
H^{1}(X/G,{\cal O}_{X/G}(mK_{X/G,orb}+L)\otimes{\cal I}(h^{m})\otimes{\cal I}_{(X/G)_{1}}) = 0
\]
holds. 
This completes the proof of Lemma 3.3.
\vspace{5mm}{\bf Q.E.D.} \\ 

Let $\tau$ be a general section in 
$H^{0}(X/G,{\cal O}_{X/G}(L))$.

Let $m_{1}$ be a sufficiently large positive integer divisible by $\mid G\mid$ 
and let $\sigma_{1}^{\prime}$
be a general element of 
\[
H^{0}((X/G)_{1},{\cal O}_{(X/G)_{1}}(m_{1}K_{X/G,orb})\otimes{\cal I}(h^{m_{1}})\otimes
{\cal M}_{x_{1},x_{2}}^{\lceil\sqrt[n_{1}]{\mu_{1}}(1-\varepsilon )\frac{m_{1}}
{\sqrt[n_{1}]{2}}\rceil}),
\]
where $x_{1},x_{2}\in (X/G)_{1}$ are distinct nonsingular points on $(X/G)_{1}$.

By Lemma 3.2, we may assume that $\sigma_{1}^{\prime}$ is nonzero.
Then by Lemma 3.3 we see that   
\[
\sigma_{1}^{\prime}\otimes\tau\in
H^{0}((X/G)_{1},{\cal O}_{(X/G)_{1}}(m_{1}K_{X/G,orb}+L)\otimes{\cal I}(h^{m_{1}})\otimes
{\cal M}_{x_{1},x_{2}}^{\lceil\sqrt[n_{1}]{\mu_{1}}(1-\varepsilon )\frac{m_{1}}
{\sqrt[n_{1}]{2}}\rceil})
\]
extends to a section
\[
\sigma_{1}\in H^{0}(X/G,{\cal O}_{X/G}((m+\ell )K_{X/G,orb})
\otimes{\cal I}(h^{m+\ell}))
\]
We may assume that  there exists a neighbourhood $U_{x,y}$ of $\{ x,y\}$ such that the divisor $(\sigma _{1})$  is smooth
on  $U_{x,y} - (X/G)_{1}$ by Bertini's theorem, if we take $\ell$
sufficiently large, since by Theorem 2.1, 
\[
H^{0}(X/G,{\cal O}_{X/G}(mK_{X/G,orb}+L)\otimes{\cal I}(h^{m}))
\rightarrow
\]
\[
\hspace{40mm}
H^{0}(X/G,{\cal O}_{X/G}(mK_{X/G,orb}+L)\otimes{\cal I}(h^{m}))/
{\cal O}_{X/G}(-(X/G)_{1})\cdot{\cal M}_{y})
\]
is surjective for every $y\in X/G$ and
 $m\geq 0$ divisible by $\mid G\mid$, where ${\cal O}_{X/G}(-(X/G)_{1})$
is the ideal sheaf of $(X/G)_{1}$.
We define a singular hermitian metric $h_{1}$ on $K_{X/G,orb}$ by
\[
h_{1} = \frac{1}{\mid\sigma_{1}\mid^{\frac{2}{m_{1}+\ell}}}.
\]
Let $\varepsilon_{0}$ be a sufficiently small positive number and 
let ${\cal I}_{1}(\alpha )$ be the multiplier ideal sheaf of 
$h_{0}^{\alpha_{0}-\varepsilon_{0}}\cdot h_{1}^{\alpha}$,i.e.,
\[
{\cal I}_{1}(\alpha ) = {\cal L}^{2}({\cal O}_{X/G},
h_{0}^{\alpha_{0}-\varepsilon_{0}}h_{1}^{\alpha}/
h_{X/G}^{(\alpha_{0}+\alpha-\varepsilon_{0})}).
\]
Suppose that $x,y$ are nonsingular points on $(X/G)_{1}$.
Then we set $x_{1} = x, x_{2} = y$ and define $\alpha_{1}(=\alpha_{1}(x,y))> 0$ by
\[
\alpha_{1} := \inf\{\alpha\mid ({\cal O}_{X/G}/{\cal I}_{1}(\alpha ))_{x}
\neq 0\,\mbox{and}\, ({\cal O}_{X/G}/{\cal I}_{1}(\alpha ))_{y}\neq 0\}.
\]
By Lemma 3.3 we may assume that we have taken $m_{1}$ so that  
\[
\frac{\ell}{m_{1}} \leq 
\varepsilon_{0}\frac{\sqrt[n_{1}]{\mu_{1}}}{n_{1}\sqrt[n_{1}]{2}}
\]
holds.
\begin{lemma}
\[
\alpha_{1}\leq \frac{n_{1}\sqrt[n_{1}]{2}}{\sqrt[n_{1}]{\mu_{1}}} 
+ O(\varepsilon _{0})
\]
holds.
\end{lemma}
To prove Lemma 3.4, we need the following elementary lemma.
\begin{lemma}(\cite[p.12, Lemma 6]{t})
Let $a,b$ be  positive numbers. Then
\[
\int_{0}^{1}\frac{r_{2}^{2n_{1}-1}}{(r_{1}^{2}+r_{2}^{2a})^{b}}
dr_{2}
=
r_{1}^{\frac{2n_{1}}{a}-2b}\int_{0}^{r_{1}^{-{2}{a}}}
\frac{r_{3}^{2n_{1}-1}}{(1 + r_{3}^{2a})^{b}}dr_{3}
\]
holds, where 
\[
r_{3} = r_{2}/r_{1}^{1/a}.
\]
\end{lemma}
{\em Proof of Lemma 3.3.}
Let $(z_{1},\ldots ,z_{n})$ be a local coordinate on a 
neighbourhood $U$ of $x$ in $X/G$ such that 
\[
U \cap (X/G)_{1} = 
\{ q\in U\mid z_{n_{1}+1}(q) =\cdots = z_{n}(q)=0\} .
\] 
We set $r_{1} = (\sum_{i=n_{1}+1}^{n}\mid z_{1}\mid^{2})^{1/2}$ and 
$r_{2} = (\sum_{i=1}^{n_{1}}\mid z_{i}\mid^{2})^{1/2}$.
Then there exists a positive constant $C$ such that 
\[
\parallel\sigma_{1}\parallel^{2}\leq 
C(r_{1}^{2}+r_{2}^{2\lceil\sqrt[n_{1}]{\mu_{1}}(1-\varepsilon )\frac{m_{1}}
{\sqrt[n_{1}]{2}}\rceil})
\]
holds on a neighbourhood of $x$, 
where $\parallel\,\,\,\,\parallel$ denotes the norm with 
respect to $h_{X/G}^{m_{1}+\ell}$.
We note that there exists a positive integer $M$ such that 
\[
\parallel\sigma\parallel^{-2} = O(r_{1}^{-M})
\]
holds on a neighbourhood of the generic point of $U\cap (X/G)_{1}$,
where $\parallel\,\,\,\,\parallel$ denotes the norm with respect to 
$h_{X/G}^{m_{0}}$. 
Then by Lemma 3.5, we have the inequality 
\[
\alpha_{1} \leq (\frac{m_{1}+\ell}{m_{1}})\frac{n_{1}\sqrt[n_{1}]{2}}{\sqrt[n_{1}]{\mu_{1}}} 
+ O(\varepsilon _{0})
\] 
holds. 
By using the fact that 
\[
\frac{\ell}{m_{1}} \leq 
\varepsilon_{0}\frac{\sqrt[n_{1}]{\mu_{1}}}{n_{1}\sqrt[n_{1}]{2}}
\]
we obtain that 
\[
\alpha_{1}\leq \frac{n_{1}\sqrt[n_{1}]{2}}{\sqrt[n_{1}]{\mu_{1}}} 
+ O(\varepsilon _{0})
\]
holds.
Q.E.D. \vspace{10mm} \\
If $x$ or $y$ is a singular point on $(X/G)_{1}$, we need the following lemma.
\begin{lemma}
Let $\varphi$ be a plurisubharmonic function on $\Delta^{n}\times{\Delta}$.
Let $\varphi_{t}(t\in\Delta )$ be the restriction of $\varphi$ on
$\Delta^{n}\times\{ t\}$.
Assume that $e^{-\varphi_{t}}$ does not belong to $L^{1}_{loc}(\Delta^{n},O)$
for every $t\in \Delta^{*}$.

Then $e^{-\varphi_{0}}$ is not locally integrable at $O\in\Delta^{n}$.
\end{lemma}
Lemma 3.6 is an immediate consequence of \cite{o-t}.
Using Lemma 3.6 and Lemma 3.5, we see that Lemma 3.4 holds
by letting $x_{1}\rightarrow x$ and $x_{2}\rightarrow y$.

\vspace{10mm}

For the next we consider Case 1.2 and Case 2.2.  
We note that in Case 2.2 by modifying $\sigma$ a little bit 
, if necessary we may assume that
$({\cal O}_{X/G}/{\cal I}(\alpha_{0}-\varepsilon ))_{y}\neq 0$ 
and $({\cal O}_{X/G}/{\cal I}(\alpha_{0}-\varepsilon^{\prime}))_{x} = 0$ hold
for a sufficiently small positive number $\varepsilon^{\prime}$. 
For example it is sufficient to replace $\sigma$ by 
the following $\sigma^{\prime}$ constructed below.

Let $X/G^{\prime}_{1}$ be a subvariety which defines a branch of 
\[
\mbox{Spec}({\cal O}_{X/G}/{\cal I}(\alpha +\delta))
\]
at $y$.  By the assumption (changing $(X/G)_{1}$, if necessary) we may assume that $(X/G)_{1}^{\prime}$ does not 
contain $x$.  
Let $m^{\prime}$ be a sufficiently large positive integer divisible by $\mid G\mid$ such that 
$m^{\prime}/m_{0}$ is sufficiently small (we can take $m_{0}$ 
arbitrary large). 

Let $\tau_{y}$ be a general element of 
\[
H^{0}(X/G,{\cal O}_{X/G}(m^{\prime}K_{X/G,orb})\otimes 
{\cal I}_{(X/G)_{1}^{\prime}}), 
\]
where ${\cal I}_{(X/G)_{1}^{\prime}}$ is the ideal sheaf of 
$(X/G)_{1}^{\prime}$. 
If we take $m^{\prime}$ sufficiently large,
 $\tau_{y}$ is not identically zero. 
We set 
\[
\sigma^{\prime} = \sigma\cdot\tau_{y}. 
\] 
Then we see that the new singular hermitian metric $h_{0}^{\prime}$
defined by $\sigma^{\prime}$ satisfies the desired property.

In these cases, instead of Lemma 3.2, we use the following simpler lemma.

\begin{lemma} Let $\varepsilon$ be a sufficiently small positive number and let $x_{1}$ be a smooth point on $(X/G)_{1}$. 
Then for a sufficiently large $m >1$ divisible by $\mid G\mid$,
\[
H^{0}((X/G)_{1},{\cal O}_{(X/G)_{1}}(mK_{X/G,orb})\otimes{\cal I}(h^{m})\otimes
{\cal M}_{x_{1}}^{\lceil\sqrt[n_{1}]{\mu_{1}}(1-\varepsilon )m\rceil})\neq 0
\]
holds.
\end{lemma}

Then taking a general $\sigma_{1}^{\prime}$ in
\[
H^{0}((X/G)_{1},{\cal O}_{(X/G)_{1}}(m_{1}K_{X/G,orb})\otimes{\cal I}(h^{m_{1}})\otimes
{\cal M}_{x_{1}}^{\lceil\sqrt[n_{1}]{\mu_{1}}(1-\varepsilon )m_{1}
\rceil}),
\]
for a sufficiently large $m_{1}$.
As in Case 1.1 and Case 2.1 we obtain a proper subvariety
$(X/G)_{2}$ in $(X/G)_{1}$ also in this case.

Inductively for distinct points $x,y\in X/G^{\circ}$, we construct a strictly decreasing
sequence of subvarieties
\[
X/G = (X/G)_{0}(x,y)\supset (X/G)_{1}(x,y)\supset \cdots 
\]
\[
\hspace{40mm} \supset (X/G)_{r}(x,y)\supset (X/G)_{r+1}(x,y) = \{ x\}\,\mbox{or}\, \{ x,y\} ,
\]
where $R_{y}$ (or $R_{x}$) is a subvariety such that $x$ deos not
belong to $R_{y}$ and $y$ belongs to $R_{y}$.
and invariants :
\[
\alpha_{0}(x,y) ,\alpha_{1}(x,y),\ldots ,\alpha_{r}(x,y),
\]
\[
\mu_{0},\mu_{1}(x,y),\ldots ,\mu_{r}(x,y)
\]
and
\[
n >  n_{1}> \cdots > n_{r}.
\]
By Nadel's vanishing theorem (Theorem 2.1) we have the following lemma.
\begin{lemma} 
Let $x,y$ be two distinct points on $X/G^{\circ}$. 
Then for every $m\geq \lceil\sum_{i=0}^{r}\alpha_{i}(x,y)\rceil +1$,
$\Phi_{\mid mK_{X/G,orb}\mid}$ separates $x$ and $y$.
\end{lemma}
{\em Proof}. 
For simplicity let us denote $\alpha_{i}(x,y)$ by $\alpha_{i}$. 
Let us define the singular hermitian metric $h_{x,y}$ of the {\bf Q}-line bundle $(m-1)K_{X/G,orb}$ defined by  
\[
h_{x,y} = (\prod_{i=0}^{r-1}h_{i}^{\alpha_{i}-\varepsilon_{i}})\cdot
 h_{r}^{\alpha_{r}+\varepsilon_{r}}h^{(m-1-(\sum_{i=0}^{r-1}(\alpha_{i}-\varepsilon_{i}))- (\alpha_{r}+\varepsilon_{r})-\ell\delta_{L})}\cdot h_{L}^{\delta_{L}},
\]
where $h_{L}$ is a $C^{\infty}$-hermitian metric on $L$ with strictly positive curvature and $\delta_{L}$ be a sufficiently small positive number. 
Then we see that  ${\cal I}(h_{x,y})$ defines a subscheme of 
$X/G$ with isolated support around $x$ or $y$ by the definition of 
the invariants $\{\alpha_{i}\}$'s. 
By the construction the curvature current $\Theta_{h_{x,y}}$ is strictly positive on $X/G$. 
Then by Nadel's vanishing theorem (Theorem 2.1) we see that 
\[
H^{1}(X/G,{\cal O}_{X/G}(K_{X/G} + \lceil (m-1)K_{X/G,orb}\rceil )\otimes {\cal I}(h_{x,y})) = 0.
\]
Hence 
\[
H^{0}(X/G,{\cal O}_{X/G}(K_{X/G}+(m-1)\lceil K_{X/G,orb}\rceil ))
\]
separates $x$ and $y$. 
We note that 
\[
H^{0}(X/G,{\cal O}_{X/G}(K_{X/G}+(m-1)\lceil K_{X/G,orb}\rceil ))
\]
is a subspace of 
\[
H^{0}(X,{\cal O}_{X}(mK_{X}))^{G}
\]
by the definition of $K_{X/G,orb}$. 
This implies that $\Phi_{\mid [mK_{X/G,orb}]\mid}$ separates 
$x$ and $y$.  {\bf  Q.E.D.} \vspace{10mm} \\

We note that for a fixed $x$, $\sum_{i=0}^{r}\alpha_{i}(x,y)$ depends on $y$.
We set
\[
\alpha (x) = \sup_{y\in U_{0}}\sum_{i=0}^{r}\alpha_{i}
\]
and let 
\[
X/G = (X/G)_{0}\supset (X/G)_{1}\supset (X/G)_{2}\supset\cdots 
\]
\[
\hspace{40mm} (X/G)_{r} \supset
(X/G)_{r+1} = \{ x\}\,\mbox{or}\,\{ x,y\}
\]
be the stratification which attains $\alpha (x)$.
In this case we call it the maximal stratification at $x$.
We see that there exists a nonempty open subset $U$ 
in countable Zariski topology 
of $X/G$ such that on $U$ the function $\alpha (x)$ is constant
and there exists an irreducible family of stratification 
which attains $\alpha (x)$ for every $x\in U$.

In fact this can be verified as follows.
We note that the cardinarity of
\[
\{ (X/G)_{i}(x,y)\mid  x,y \in X/G, x\neq y (i=0,1,\ldots )\}
\]
is uncontably many, while the cardinarity of the irreducible components
of Hilbert scheme of $X/G$ is countably many.
We see that for fixed $i$ and very general $x$, $\{ (X/G)_{i}(x,y)\}$ should form a family on $X/G$.  Similary we see that for very general $x$, we may assume that
the maximal stratification $\{ (X/G)_{i}(x)\}$ forms a family.  
 This implies the existence of $U$.   

And we may also assume that the corresponding invariants $\{\alpha_{0},
\ldots ,\alpha_{r}\}$, $\{\mu_{0},\ldots ,\mu_{r}\}$,
$\{ n = n_{0}\ldots ,n_{r}\}$ are constant on $U$.
Hereafter we denote these invariants again  by the same notations for simplicity.
The proof of the following lemma is parallel to that of Lemma 3.4.
\begin{lemma}
\[
\alpha_{i}\leq \frac{n_{i}\sqrt[n_{i}]{2}}{\sqrt[n_{i}]{\mu_{i}}} + O(\varepsilon_{i-1})
\]
hold for $1\leq i\leq r$.
\end{lemma}

\begin{proposition}
For every 
\[
m > \lceil\sum_{i=0}^{r}\alpha_{i}\rceil + 1
\]
$\mid [mK_{X/G,orb}]\mid$ gives a birational rational map from $X/G$ into 
a projective space.
\end{proposition}

\begin{lemma} If $\Phi_{m}\mid_{(X/G)_{i}}$ is birational rational map
onto its image, then
\[
\deg \Phi_{m}((X/G)_{i})\leq m^{n_{i}}\mu_{i}
\]
holds.
\end{lemma}
{\em Proof}.
Let $p : \tilde{X/G}\longrightarrow X/G$ be the resolution of 
the base locus of $\mid mK_{X/G,orb}\mid$ and let 
\[
p^{*}\mid [mK_{X/G,orb}]\mid = \mid P_{m}\mid + F_{m}
\]
be the decomposition into the free part $\mid P_{m}\mid$ 
and the fixed component $F_{m}$. 
Let $p_{i} : \tilde{X/G}_{i}\longrightarrow (X/G)_{i}$ be the resolution
of the base locus of $\Phi_{\mid mK_{X/G,orb}\mid}\mid_{(X/G)_{i}}$ 
obtained by the restriction of $p$ on $p^{-1}((X/G)_{i})$. 
Let 
\[
p_{i}^{*}(\mid mK_{X/G,orb}\mid_{(X/G)_{i}}) = \mid P_{m,i}\mid + F_{m,i}
\]
be the decomposition into the free part $\mid P_{m,i}\mid$ and 
the fixed part $F_{m,i}$.
We have
\[
\deg \Phi_{\mid [mK_{X/G,orb}]\mid}((X/G)_{i}) = P_{m,i}^{n_{i}}
\]
holds.
Then by the ring structure of $R(X/G,K_{X/G,orb})$, we have that
there exists a natural injection
\[
H^{0}(X/G,{\cal O}_{X/G}(\nu P_{m}))\rightarrow 
H^{0}(X/G,{\cal O}_{X/G}([m\nu K_{X/G,orb}])\otimes{\cal I}(h^{m\nu}))
\]
for every $\nu\geq 1$.
Hence there exists a natural morphism
\[
H^{0}((X/G)_{i},{\cal O}_{(X/G)_{i}}(\nu P_{m,i}))
\rightarrow 
H^{0}((X/G)_{i},{\cal O}_{(X/G)_{i}}([m\nu K_{X/G,orb}])\otimes{\cal I}(h^{m\nu}))
\]
for every $\nu\geq 1$. 
This morphism is clearly injective. 
This implies that 
\[
\mu_{i} \geq  m^{-n_{i}}\mu ((X/G)_{i},P_{m,i})
\]
holds. 
Since $P_{m,i}$ is nef and big on $(X/G)_{i}$ we see that 
\[
\mu ((X/G)_{i},P_{m,i}) = P_{m,i}^{n_{i}}
\]
holds.
Hence
\[
\mu_{i}\geq m^{-n_{i}}P_{m,i}^{n_{i}}
\]
holds.  This implies that
\[
\deg \Phi_{\mid mK_{X/G,orb}\mid}((X/G)_{i})\leq \mu_{i}m^{n_{i}}
\]
holds.
{\bf Q.E.D.}

\section{Proof of Theorem 1.1}
To prove Theorem 1.1 we use the following subadjunction formula. 
\begin{theorem}(\cite{ka})
Let $X/G$ be a normal projective variety.
Let $D^{\circ}$ and $D$ be effective {\bf Q}-divisor on $X$ such that 
$D^{\circ} < D$, $(X,D^{\circ})$ is logterminal and 
$(X,D)$ is logcanonical. 
Let $W$ be a minimal center of logcanonical singularities for $(X,D)$. 
Let $H$ be an ample Cartier divisor on $X$ and $\epsilon$ a positive rational number.
Then there exists an effective {\bf Q}-divisor $D_{W}$ on $D$ such that 
\[
(K_{X}+D+\epsilon H)\mid_{W}\sim_{\mbox{\bf Q}}K_{W}+D_{W}
\]
and $(W,D_{W})$ is logterminal. 
In particular $W$ has only rational singularities.
\end{theorem} 
Let us start the proof of Theorem 1.1. 
We prove Theorem 1.1 by induction on $n = \dim X$. 
Suppose that Theorem 1.1 holds for varities of general type of dimension $< n$.
Then there exists a positive constant $C(m) (m < n)$ depending only on $m$ such that for every smooth projective varietiey $Y$ of general type of dimension $m$
\[
\mu (Y,K_{Y})/\sharp\mbox{Aut}(Y)  \geq C(m)
\]
holds. 
Let $X$ be a smooth projective variety of general type as in Section 3. 
We use the same notations as in Section 3. 
Let $x,y$ be distinct points on $(X/G)^{\circ}$ and 
let 
\[
X/G = (X/G)_{0} \supset (X/G)_{1}\supset\cdots (X/G)_{r}\supset (X/G)_{r+1} = \{ x\} \mbox{or}\,\,\{ x,y\}
\]
be the stratification constructed as in Section 3 and let 
\[
\mu_{0},\ldots , \mu_{r}
\]
\[
n_{1},\ldots ,n_{r}
\]
be the invariants as in Section 3. 
Let 
\[
X= X_{0} \supset X_{1}\supset\cdots X_{r}\supset X_{r+1}
\]
be the corresponding stratification of $X$. 
If we take $x,y$ general, $X_{i} (0\leq i\leq r)$ are projective varieties of 
general type. 
Let 
\[
X_{can} := \mbox{Proj}\, R(X,K_{X})
\]
be the canonical model of $X$.
 
We have the corresponding stratification 
\[
X_{can} = X_{0,can} \supset X_{1,can}\supset\cdots X_{r,can}\supset X_{r+1,can}
\]
on $X_{can}$ (here we note that $X_{i,can}$ does not denote the canonical model of $X_{i}$ for $i\geq 1$). 

Then we see that 
\[
\mu_{i} = \frac{1}{\mid G\mid}\mu (X_{i},K_{X}) = \frac{1}{\mid G\mid}(K_{X_{can}})^{n_{i}}\cdot X_{i,can}
\]
holds. 
Let $H$ be an ample divisor on $X$. 
By the subadjunction formula, we see that for every positive rational number $\epsilon$ 
\[
K_{X_{i,can}} <_{\mbox{\bf Q}} (1+\sum_{j=0}^{i-1}\alpha_{j})K_{X_{can}}+\epsilon H
\]
holds, where $<_{\mbox{\bf Q}}$ means that the righthandside minus the lefthandside is {\bf Q}-linear equivalent to an effective divisor and 
$K_{X_{i},can}$ denotes the pushforward of the canonical divisor of a nonsingular model of $K_{X_{i,can}}$.  
This can be verified as follows. 
Let 
\[
\pi : X \longrightarrow X/G
\]
be the natural morphism. 
Let $D_{i}$ be the divisor on $X$ which corresponds to 
the singular hermitian metric 
\[
\pi^{*}(h_{0}^{\alpha_{0}-\varepsilon_{0}}\cdots h_{i-1}^{\alpha_{i-1}-\varepsilon_{i-1}}\cdot h_{i}^{\alpha_{i}}).
\]
$D_{i}$ is a positive linear combinations of $\{\pi^{*}(\sigma_{0}),\ldots ,(\sigma_{j})\}$ by the constructions of $h_{0},\ldots , h_{i}$. 
Also we may assume that $D_{i}$ is a ${\bf Q}$-divisor by perturbations of 
$\varepsilon_{0},\ldots ,\varepsilon_{i-1}$. 
 $X_{i,can}$ may not be the minimal center of $(X,D_{i})$ and $(X,D_{i})$  may not be logcanonical. 
But if we take a suitable modification 
\[
\pi_{i} : Y_{i} \longrightarrow X_{i,can},
\]
we may assume that there exists an effective {\bf Q}-divisor $E_{i}$ such that 
\begin{enumerate}
\item 
$\pi_{i}^{*}D_{i}-E_{i}$ is effective,
\item 
$(Y_{i},\pi_{i}^{*}D_{i}- E_{i})$ is logcanonical and the proper transform of 
$X_{i,can}$ is the minimal center of $(Y_{i},\pi_{i}^{*}D_{i}- E_{i})$. 
\end{enumerate}
Then by Theorem 4.1, we have that for every positive rational number $\epsilon$
\[
K_{X_{i,can}} <_{\mbox{\bf Q}} (1+\sum_{j=0}^{i-1}\alpha_{j})K_{X_{can}}+\epsilon H
\]
holds.
By the inductive assumption this implies that 
\[
(1 +\sum_{j=0}^{i-1}\alpha_{j})^{n_{i}}\cdot \mu_{i}\geq C(n_{i})
\]
holds. 
Since 
\[
\alpha_{i} \leq \frac{\sqrt[n_{i}]{2}n_{i}}{\sqrt[n_{i}]{\mu_{i}}}+ O(\varepsilon_{i-1})
\]
holds by Lemma 3.9, 
we see that 
\[
(*) \hspace{10mm} 
\frac{1}{\sqrt[n_{i}]{\mu_{i}}}\leq (1+\sum_{j=0}^{i-1}\frac{\sqrt[n_{j}]{2}n_{j}}{\sqrt[n_{j}]{\mu_{j}}})\cdot C(n_{i})^{-1}
\]
holds for every $i \geq 1$.
Inductively we see that if $\mu_{0}\leq 1$ holds,
\[
\frac{1}{\sqrt[n_{i}]{\mu_{i}}}\leq \frac{1}{\sqrt[n]{\mu_{0}}}C(C(1),\ldots ,C(n-1))
\]
holds where $C(C(1),\ldots C(n-1))$ is a positive constant depending only on 
$C(1),\ldots ,C(n-1)$.
Hence if $\mu_{0} < 1$ holds then we see that 
\[
\deg \Phi_{\mid (1+\lceil\sum_{i=0}^{r}\alpha_{i}\rceil)K_{X/G,orb}\mid}(X)
\leq C(C(1),\ldots ,C(n-1))^{n}
\]
holds.
This implies that $X/G$ is birationally bounded, if 
\[
\mu_{0} (= \frac{1}{\mid G\mid}\mu (X,K_{X})) \leq 1
\]
holds.
We set 
\[
\alpha := \lceil\sum_{i=0}^{r}\alpha_{i}+1\rceil .
\]
Then using Lemma 3.10, we have the following lemma.
\begin{lemma} If $\mu_{0}\leq 1$ holds, then  
there exists a positive constant $A(n)$ depending only on $n$ 
such that 
\[
1\leq \alpha^{n}\mu_{0} \leq A(n)
\]
holds. 
\end{lemma}
Let 
\[
\mid \alpha K_{X}\mid^{G} = \mid P\mid + F
\]
be the decomposition of $\mid \alpha K_{X}\mid^{G}$ 
into the movable part $\mid P\mid$ and the fixed 
component $F$. 
Taking a suitable successive $G$-equivariant blowing ups, 
we may assume that $\mid P\mid$ is base point free. 
And also we may assume that 
the canonical birational map 
\[
f : X \longrightarrow X_{can}
\]
is a morphism. 
\begin{lemma}
There exists a positive constant $c_{n}$ depending only 
on $n$ such that 
\[
f^{*}K_{X_{can}}\cdot P^{n-1}\geq c_{n}\mid G\mid
\]
holds.
In particular 
\[
\alpha^{n-1}K_{X_{can}/G,orb}^{n}\geq c_{n}
\]
holds.
\end{lemma}
{\em Proof}.
Let 
\[
f_{G} : X/G \longrightarrow X_{can}/G
\]
be the natural morphism. 
Let us write
\[
K_{X/G} = f_{G}^{*}(K_{X_{can}/G}) + \sum a_{i}E_{i}
\]
where $\{ E_{i}\}$ are irreducible exceptional divisor of $f_{G}$.
We set 
\[
Y := \Phi_{\mid \alpha K_{X}\mid^{G}}(X).
\]
and we set 
\[
\phi := \Phi_{\mid P\mid} : X \longrightarrow Y.
\]
Let 
\[
\phi_{G} : X/G \longrightarrow Y
\]
be the birational morphism induced by $\phi$. 
Then 
\[
 f^{*}K_{X_{can}}\cdot P^{n-1}
= \phi_{*}f^{*}K_{X_{can}}\cdot H^{n-1}
\]
holds, where $H$ denotes the hyperplane section of $Y$. 
Also 
\[
\phi_{*}f^{*}K_{X_{can}}\cdot H^{n-1} 
= \mid G\mid\cdot (\phi_{G})_{*}f_{G}^{*}K_{X_{can}/G,orb}\cdot H^{n-1}
\]
holds. 
On the other hand
\begin{eqnarray*}
(*) \hspace{10mm} (\phi_{G})_{*}f_{G}^{*}K_{X_{can}/G}\cdot H^{n-1}
 & =  & (\phi_{G})_{*}(K_{X/G}-\sum a_{i}E_{i})\cdot H^{n-1} \\
& = & K_{Y}\cdot H^{n-1} - \sum_{i}a_{i}(\phi_{G})_{*}E_{i}\cdot H^{n-1}
\end{eqnarray*}
holds, where $K_{Y}$ denotes the pushforward of the canonical divisor of 
the normalization of $Y$ to $Y$.
We note that $K_{Y}\cdot H^{n-1}(= K_{X/G}\cdot P^{n-1}$)is an integer. 
Since $E_{i}$'s appear as fixed components of $\mid  [\alpha K_{X_{can}/G,orb}]\mid^{G}$, 
we see that 
\[
\sum_{i}(\phi_{G})_{*}E_{i}\cdot H^{n-1} \leq  \alpha^{n}\mu_{0} \leq C(n)
\]
hold. 
Hence $\sum_{i}(\phi_{G})_{*}E_{i}$ is bounded.
 
Since $\sum_{i}(\phi_{G})_{*}E_{i}$ is an exceptional divisor of 
the birational rational map 
\[
f_{G}\circ\phi_{G}^{-1} : Y -\cdots\rightarrow X_{can}/G,
\]
$\{ a_{i}\}$ is  of finitely many possibilities. 
Hence there exists a positive constant $K_{n}$ depending only on $n$ such that
\[
(\sharp )\hspace{10mm} (\phi_{G})_{*}f_{G}^{*}(K_{X_{can}/G})\cdot H^{n-1} \geq -K_{n}
\]
holds.
Let $\{ D_{j}\}$ be the irreducible divisors  such that 
\[
K_{X_{can}/G,orb} = K_{X_{can}/G}+\sum_{j}\frac{m_{j}-1}{m_{j}}D_{j}
\]
for some positive integers $\{ m_{j}\}$.
Then we see that 
\begin{eqnarray*}
(\flat ) \,\,\,(f_{G}^{*}K_{X_{can}/G,orb})\cdot \phi_{G}^{*}H^{n-1}
& = & f_{G}^{*}K_{X_{can}/G}\cdot  \phi_{G}^{*}H^{n-1} 
+ \sum_{j}\frac{m_{j}-1}{m_{j}}f_{G}^{*}D_{j}\cdot \phi_{G}^{*}H^{n-1} \\
 & \leq  & \alpha^{n}\mu_{0} \\  & \leq  & A(n)
\end{eqnarray*}
hold. 
By $(\sharp )$ this implies that $\sum_{j}(\phi_{G})_{*}f_{G}^{*}D_{j}$ 
is bounded and 
\[
\sharp \{ j \mid  (\phi_{G})_{*}f_{G}^{*}D_{j}\neq 0 \}
\]
is uniformly bounded by a positive integer, say $N$ depending only on $n$. 

\begin{lemma}
Let $N$ and $B$ are fixed positive integers. 
Then
\[
\{ \{ -\sum_{j=1}^{N}\frac{b_{j}}{a_{j}}\} \mid, a_{j}, b_{j} \,\,\mbox{are integers such that} \,\, b_{j}\leq B \} -
\{ 0\}
\]
is bounded below by a positive constant, where for a rational number $c$ $\{ c\}$
denotes the fractional part  of $c$. i.e. 
\[
\{ c\} := c - [c].
\]
\end{lemma}
{\em Proof}. 
Suppose not. 
Then there exists a sequence of positive integers 
\[
\{ a_{j,k}\} ,\{ b_{j,k}\} 1\leq j\leq N, k= 1,2,\ldots 
\]
such that 
\[
b_{j,k} \leq B, 
\]
\[
\{ -\sum_{j=1}^{N}\frac{b_{j,k}}{a_{j,k}}\} \neq 0,
\]
\[
\lim_{k\rightarrow\infty}\frac{b_{j,k}}{a_{j,k}}
\]
exists for every $j$ and
\[
\lim_{k\rightarrow \infty} \{ -\sum_{j=1}^{N}\frac{b_{j,k}}{a_{j,k}}\}  = 0
\]
hold. 
We note that if 
\[
\lim_{k\rightarrow\infty}\frac{b_{j,k}}{a_{j,k}} \neq  0
\]
then by the boundedness of $b_{j,k}$ the sequence is constant for 
every sufficiently large $k$ 
and if 
\[
\lim_{k\rightarrow\infty}\frac{b_{j,k}}{a_{j,k}} =  0
\]
then $a_{j,k}$ tends to infinity $k$ goes to infinity. 
Since 
\[
\lim_{k\rightarrow \infty} \{ -\sum_{j=1}^{N}\frac{b_{j,k}}{a_{j,k}}\}  = 0
\]
holds, there is no $j$ such that 
\[
\lim_{k\rightarrow\infty}\frac{b_{j,k}}{a_{j,k}} =  0
\]
holds. 
Hence by the above observation we see that for every $j$ 
the sequence $\{ b_{j,k}/a_{j,k}\}_{k=1}^{\infty}$ 
is constant  
for every sufficiently large $k$ and $j$.
This contradicts to the fact that 
\[
\{ -\sum_{j=1}^{N}\frac{b_{j,k}}{a_{j,k}}\} \neq 0
\]
holds for every $k$.
This completes the proof of Lemma 4.3 . {\bf Q.E.D.} \vspace{5mm} \\

We note that by $(*)$, the finiteness properties of $\{ a_{i}\}$ 
and the boundedness of $\sum_{i}(\phi_{G})_{*}E_{i}$, we see that 
the rational number 
$f_{G}^{*}K_{X_{can}/G}\cdot H^{n-1}$ is of finitely many possibilities. 
By $(\flat )$, the boundedness of  $\sum_{j}(\phi_{G})_{*}f_{G}^{*}D_{j}$ and 
Lemma 4.3, we see that 
there exists a positive constant $c_{n}$ depending only 
on $n$ such that 
\[
f^{*}K_{X_{can}}\cdot P^{n-1}\geq c_{n}\mid G\mid
\]
holds.
Since $R(X_{can}/G,K_{X_{can}/G,orb})$ is a ring, 
\[
\alpha^{n-1}K_{X_{can}/G,orb}^{n}\geq c_{n}
\]
holds.
This completes the proof of Lemma 4.1. {\bf Q.E.D.} \vspace{5mm} \\
By Lemma 4.1 and Lemma 4.2 we see that 
\[
\alpha \leq \frac{A(n)}{c_{n}}
\]
holds. 
By Lemma 4.1, we see that 
\[
\mu_{0} \geq \frac{1}{\alpha^{n}}
\]
holds.
Hence we have that 
\[
\mu_{0} \geq (\frac{c_{n}}{A(n)})^{n}
\]
holds.
This completes the proof of Theorem 1.2.
Since 
\[
\mu_{0} = \frac{1}{\mid G\mid}\mu (X,K_{X})
\]
holds, we have that
\[
\mid G\mid \leq (\frac{A(n)}{c_{n}})^{n}\mu (X,K_{X})
\]
holds.
This completes the proof of Theorem 1.1.

Author's address\\
Hajime Tsuji\\
Department of Mathematics\\
Tokyo Institute of Technology\\
2-12-1 Ohokayama, Megro 152-8551\\
Japan \\
e-mail address: tsuji@math.titech.ac.jp

\end{document}